\documentclass{article}
\usepackage{latexsym}
\usepackage{amsthm}  
\usepackage[all]{xy}
\usepackage{stmaryrd}

\newcommand{\mbb}[1]{\mathbf{#1}}

\newcommand{\msc}[1]{\mathcal{#1}}
\newcommand{\mrm}[1]{\mathrm{#1}}
\newcommand{\C}{\mbb{C}}

\newcommand{\Q}{\mbb{Q}}

\newcommand{\Pic}{\mrm{Pic~}}

\newcommand{\tens}{\varotimes}
\newcommand{\dirsum}{\varoplus}

\newcommand{\sh}[1]{{\msc{#1}}}
\newcommand{\fn}[3]{{#1:#2\rightarrow #3}}

\newcommand{\shhom}{{\msc{H}\!\!\mathit{om}}}

\newcommand{\Hom}{\mrm{Hom}}

\newcommand{\Ext}{\mrm{Ext}}
\newcommand{\EExt}{\mrm{\mbb{Ext}}}
\renewcommand{\L}{\mbb{L}}
\newcommand{\Spec}{\mrm{Spec~}}

\newcommand{\Der}{\mrm{Der}}
\newcommand{\RHom}{{\mrm{RHom}}}
\newcommand{\shRHom}{{R\shhom}}

\theoremstyle{definition}
\newtheorem{defn}{Definition}[section]

\theoremstyle{plain}
\newtheorem{thm}[defn]{Theorem}
\newtheorem{cor}[defn]{Corollary}
\newtheorem{lem}[defn]{Lemma}
\newtheorem{prop}[defn]{Proposition}

\date{}		

\begin{document}

\title{Moduli of products of curves}
\author{Michael A. van Opstall}

\maketitle

\footnotetext{University of Washington, {\tt opstall@math.washington.edu}.}

\begin{abstract}
Some technical results on the deformations of varieties of general type
and on permanence of semi-log-canonical singularities are proved. These
results are applied to show that the connected component of the moduli
space of stable surfaces containing the moduli point of a product of stable
curves is the product of the moduli spaces of the curves, assuming the
curves have different genera. An application of this result shows that
even after compactifying the moduli space and fixing numerical invariants,
the moduli spaces are still very disconnected.
\end{abstract}

The main result of this article is the construction of several connected,
irreducible components of the moduli space of stable surfaces. These components
parameterize products of stable curves. They are constructed from the
corresponding moduli spaces of stable curves. The essential results are 
that products of stable curves are stable surfaces, and all infinitesimal
deformations of a product of stable varieties (of any dimension) come from
the deformations of the factors. 
In particular, this gives a proof that even after fixing certain topological
invariants, the resulting moduli space may have arbitrarily many components.

These examples also show that these components of the moduli space of smooth
minimal surfaces of general type are not joined in the moduli space of
stable surfaces.

All schemes are defined over the field
$\C$. A {\em variety} will be a connected, reduced, and separated scheme of 
finite type, not
necessarily assumed irreducible. A {\em family} will be a flat morphism of 
varieties. The deformation theory results used can be found in 
\cite{vis:dlci};
 the base space of a miniversal deformation will be called the {\em 
Kuranishi space}, following the convenient terminology from analytic geometry.

I thank S\'andor Kov\'acs and E. Lee Stout for useful conversations. Thanks
also to an anonymous reader who noted some simplifications from an earlier 
version.

\section{Semi-log-canonical singularities and moduli of stable surfaces}

First a preliminary remark: if $X$ is a $S_2$ variety which is 
Gorenstein in codimension 1, then the extension of the dualizing sheaf
of the Gorenstein locus (which is locally free) is a reflexive sheaf on
$X$ which corresponds to a Weil divisor $K_X$. $X$ is called $\Q$-Gorenstein
if some multiple of $K_X$ is Cartier. Using these definitions, one may
define the class of singularities to be studied.

\begin{defn}
A variety $X$ is said to have {\em semi-log-canonical} (slc) singularities
if
\begin{enumerate}
\item $X$ is $\Q$-Gorenstein;
\item $X$ is $S_2$;
\item $X$ has at worst normal crossing singularities is codimension 1;
\item there exists a good desingularization $\fn{f}{Y}{X}$ such that in the
formula
\[
K_Y\equiv f^*K_X+\sum a_iE_i
\]
all of the $a_i$ are positive.
\end{enumerate}
\end{defn}

The moduli space of stable surfaces with fixed Hilbert polynomial is 
difficult to define. In fact, there are many different definitions which
make sense but lead to nonisomorphic moduli spaces. The original definition 
(as well as the definition of slc singularities) appeared in \cite{ksb:3f}.
Later Koll\'ar amended the conditions that a family of stable varieties should
satisfy in \cite{k:pcm}.

In the case of Gorenstein varieties, these subtleties do not occur, and the 
moduli functors do not differ. However, it is only in special
cases that the moduli space of minimal surfaces of general type can be
compactified by adding stable surfaces with only Gorenstein singularities to
the moduli problem. Since this suffices for the purposes of this paper, only
this special case is considered.

\begin{defn}
The moduli functor of stable Gorenstein surfaces is a functor from 
schemes to sets which assigns to a scheme $B$ the set of isomorphism
classes of flat, proper morphisms $X\rightarrow B$ whose fibers are 
Gorenstein schemes with slc singularities and whose relative dualizing
sheaf $\omega_{X/B}$ is ample. 
\end{defn}

This article considers a smaller functor. Let $M_{g_1,g_2}$ be the functor
which assigns to $B$ the set of isomorphism classes of flat proper morphisms
$X\rightarrow B$ whose fibers are products of stable curves of genera $g_1$ 
and $g_2$. The results of
this article will show that this functor is coarsely representable by a 
connected and projective variety and that it is an open and closed subfunctor 
of the moduli functor of stable Gorenstein surfaces.

\section{Deformations of products}

In this section, some general deformation-theoretic results are proved about
products of varieties. These results are formal and primarily homological.
The goal is to show that under some conditions on singularities, the small
deformations of a product of varieties are obtained by deforming the factors. 

Let $X=Y_1\times Y_2$ be a variety which is the product of two local complete
intersection varieties $Y_1$ and $Y_2$ of general type; let $\pi_i$ denote the
projection map to $Y_i$. This notation will be fixed throughout this section.

The following ``rigidity lemma'' will be useful:

\begin{lem} \label{riglem}
If $\fn{h}{X_1\times X_2}{B_1\times B_2}$ is a surjective morphism of
products of stable varieties, then after possibly renumbering, 
$h$ can be written as the product of maps $\fn{h_i}{X_i}{B_i}$, $i=1,2$.
\end{lem}

\begin{proof}
This follows from the fact that the tangent space to the scheme
$\Hom(X_i,B_j)$ at the equivalence class $[f]$ of a morphism is 
$H^0(X_i,f^*\sh{T}_{B_j})$ which vanishes due to the stability assumption. 
A morphism $h$ as in the hypothesis which is not a product would be a 
non-trivial deformation of some morphism $\fn{f}{X_i}{B_j}$.
\end{proof}

In particular, it follows that, up to renumbering, a product of curves
of general type can be written as a product of curves in a unique way. This
depends on the general type assumption, as there exist abelian surfaces 
which can be written in distinct ways as the product of elliptic curves.

The assumption that the varieties is this section are local complete
intersections implies that the space $T^1(X)$ of first-order infinitesimal
deformations of such a variety $X$ is given by $\Ext^1_X(\Omega_X,\sh{O}_X)$.
See \cite{vis:dlci} for details. Without the local complete intersection
hypothesis, the sheaf $\Omega_X$ is replaced with the cotangent complex
and $\Ext$ is taken to be the hyperext in the derived category:
$T^1(X)=\EExt^1(\L_X,\sh{O}_X)$.

\begin{thm}\label{stdef}
Every first-order deformation of $X$ is the product of a first order 
deformation of $Y_1$ with a first order deformation of $Y_2$ if $Y_1$ and
$Y_2$ are of general type.
\end{thm}

\begin{proof}
Let $\L_X$, $\L_{Y_1}$, and $\L_{Y_2}$ denote the cotangent complexes of
$X$, $Y_1$, and $Y_2$, respectively. Denote by $\EExt$ the hyperext groups.
We need to show:
\[
\EExt^1_X(\L_X,\sh{O}_X)\cong \EExt^1_{Y_1}(\L_{Y_1},\sh{O}_{Y_1})
\dirsum\EExt^1_{Y_2}(\L_{Y_2},\sh{O}_{Y_2})
\]
by \cite{ill:cc}, III.1.2.0.

By \cite{ill:cc}, II.2.2.3, 
\begin{eqnarray}
\EExt^1(\L_X,\sh{O}_X)&\cong & \EExt^1(\pi_1^*\L_{Y_1}\dirsum
\pi_2^*\L_{Y_2},\sh{O}_X) \\
&\cong & \EExt^1(\pi_1^*\L_{Y_1},\pi_1^*\sh{O}_{Y_1})\dirsum
\EExt^1(\pi_2^*\L_{Y_2},\pi_2^*\sh{O}_{Y_2}),
\end{eqnarray}

The following computation finishes the proof:
\begin{eqnarray}
\EExt^1(\pi_1^*\L_{Y_1},\pi_1^*\sh{O}_{Y_1})&\cong & 
H^1[\RHom(\pi_1^*\L_{Y_1},\pi_1^*\sh{O}_{Y_1})] \\
&\cong & H^1[R\Gamma(X,\shRHom(\pi_1^*\L_{Y_1},\pi_1^*\sh{O}_{Y_1}))] \\
&\cong & H^1[R\Gamma(X,\pi_1^*\shRHom(\L_{Y_1},\sh{O}_{Y_1}))] \\
&\cong & H^1[R\Gamma(Y_2,R{\pi_2}_*\pi_1^*\shRHom(\L_{Y_1},\sh{O}_{Y_1}))] \\
&\cong & H^1[R\Gamma(Y_2,\sh{O}_{Y_2}\tens R\Gamma(Y_1,
\shRHom(\L_{Y_1},\sh{O}_{Y_1})))] \\
&\cong & H^1[R\Gamma(Y_2,\sh{O}_{Y_2})\tens \RHom(\L_{Y_1},\sh{O}_{Y_1})] \\
&\cong &[H^0(R\Gamma(Y_2,\sh{O}_{Y_2}))\tens H^1(\RHom(\L_{Y_1},\sh{O}_{Y_1})]
\\
\nonumber & &\dirsum [H^1(R\Gamma(Y_2,\sh{O}_{Y_2}))\tens H^0(\RHom(\L_{Y_1},\sh{O}_{Y_1})]
\\
&\cong & \EExt^1(\L_{Y_1},\sh{O}_{Y_1})\dirsum [H^1(Y_2,\sh{O}_{Y_2})\tens
\Der(\sh{O}_{Y_1},\sh{O}_{Y_1})] \\
&\cong & \EExt^1(\L_{Y_1},\sh{O}_{Y_1})
\end{eqnarray}
The steps are justified as follows: the composition of derived functors
rule (\cite{hart:rd}, II.5.3) justifies steps (4), (6), and part of (8). Step 
(5) follows from the flatness of $\pi_1$ using \cite{hart:rd}, II.5.8. Step
(7) is \cite{hart:rd}, II.5.12. Step (8) follows from \cite{hart:rd}, II.5.16.
Step (9) is the K\"unneth formula. Step (10) follows from properness of
$Y_2$ and \cite{ill:cc}, II.1.2.4.3. Step (11) follows from the fact that
varieties of general type have no infinitesimal automorphisms, so
$\Der(\sh{O}_{Y_1},\sh{O}_{Y_1})$ vanishes.
\end{proof}

Note that the general type hypothesis is essential: suppose the $Y_i$ were
both smooth elliptic curves. Then one may replace all of the 
$\Ext^1(\Omega,\sh{O})$ on $X$ and the $Y_i$ with $H^1(\sh{T})$. The tangent 
sheaf of an abelian variety is trivial, so 
$h^1(Y_1,\sh{T}_{Y_1})+h^1(Y_2,\sh{T}_{Y_2})=h^1(Y_1,\sh{O}_{Y_1})+
h^1(Y_2,\sh{O}_{Y_2})=2$, but by Hodge theory, $h^1(X,\sh{T}_X)=
2h^1(X,\sh{O}_X)=\dim_\C H^1(X,\C)=4$.

\begin{cor}
The Kuranishi space of a product of finitely many stable curves is smooth.
\end{cor}

\begin{proof}
This follows from the fact that the deformations of stable curves are 
unobstructed, and from the above result shows that the
only infinitesimal deformations of the product come from the factors, and
are consequently unobstructed.
\end{proof}

\section{Products of stable varieties}

An essential advantage in using the compactified moduli space becomes
apparent when one can determine all of the stable degenerations of a class
of varieties. This section is dedicated to proving that products of smooth
curves degenerate to products of stable curves, although the result is
slightly more general. Higher-dimensional versions of such a result would
depend on a deeper study of obstructions which appear. The proof given here
also uses the normality of the moduli spaces of stable curves, which follows
from unobstructedness. 

\begin{prop}
The product of stable curves is a stable surface. More specifically, the
product of stable curves has only normal crossings and degenerate cusps as
singular points.
\end{prop}

\begin{proof}
The question is analytically local. The singularities of stable
curves are nodes. Since the product of a smooth point on a curve with a node
is a normal crossing singularity, it suffices to check that the product
of a node with itself is slc.

One must compute a semiresolution of the scheme $\Spec \C[x,y,w,z]/(xy,wz)$. 
This scheme is the affine cone over a cycle of rational curves, so blowing
up the cone point is a semiresolution with exceptional locus begin a cycle
of rational curves. Therefore the singularity at the cone point is a degenerate
cusp. All of the other singular points are plainly normal crossings.

Having checked that the singularities are slc,
the stability assertion is simply the ampleness of the canonical bundle,
which follows from the ampleness of the canonical bundles of the factors.
\end{proof}

This is a special case of the following more general result, but the proof
with coordinate rings is retained to see exactly what singularities occur
in the case of products of stable curves.

\begin{thm}
Let $Y_1$ and $Y_2$ be smoothable stable varieties. Then $Y_1\times Y_2$ is
a smoothable stable variety.
\end{thm}

\begin{proof}
The ampleness of the canonical class of $Y_1\times Y_2$ is immediate. The
smoothability is clear, since the fibered product of two smoothings will be
a smoothing, since products of rational Gorenstein singularities are
rational Gorenstein. It remains to verify that products of slc singularities
are slc. First, the conditions of $\Q$-Gorenstein, $S_2$ and normal crossings
in codimension 1 are clearly preserved under taking products.
Let $\fn{f}{X}{Y_1}$ be a desingularization. Then
write
\[
K_X=f^*K_{Y_1}+\sum a_iE_i
\] 
where the $E_i$ are exceptional. The $a_i$ are all greater than or equal
to -1 since $Y_1$ is slc. Therefore the exceptional divisors of the product
morphism $X\times Y_2\rightarrow Y_1\times Y_2$ occur with coefficient 
greater than or equal to -1. Since $X$ is smooth and $Y_2$ is slc, 
$X\times Y_2$ is slc. Therefore the discrepancies of a resolution of 
$X\times Y_2$ are all greater than or equal to -1, so $Y_1\times Y_2$ is slc,
since a resolution of $X\times Y_2$ is also a resolution of $Y_1\times Y_2$.
\end{proof}

A stronger version of this theorem which depends on minimal model 
hypotheses, and which we will not use here is in \cite{vo:thesis}. Precisely,
the total space of a flat family over a base with only slc singularities
whose special fiber has only slc singularities has only slc singularities.

\section{Main results}

The main theorems below are stated and proved in the case of the product
of two surfaces for ease of notation. However, the proofs generalize to the
product of finitely many curves. Denote by $M_{g}$ the moduli functor of
stable curves of genus $g$. This functor is known to be coarsely representable
by a projective variety.

\begin{thm}
Let $g_1,g_2\geq 2$. If $g_1\neq g_2$, then $M_{g_1,g_2}$ is isomorphic to 
$M_{g_1}\times M_{g_2}$.
\end{thm}

\begin{proof}
Taking fibered products gives a natural transformation
$M_{g_1}\times M_{g_2}\rightarrow M_{g_1,g_2}$. This natural transformation
is relatively representable. By \ref{stdef}, it is \'etale. By \ref{riglem}
it is injective on geometric points, that is $M_{g_1}(k)\times
M_{g_2}(k)=M_{g_1,g_2}(k)$ when $k$ is an algebraically closed field. The
natural transformation is proper since $M_{g_1}\times M_{g_2}$ is proper.
It follows that the functors are isomorphic and that $M_{g_1,g_2}$ is
coarsely representable.
\end{proof}

A similar argument proves:

\begin{thm}
$M_{g,g}$ is isomorphic to the symmetric square of the functor $M_g$ if
$g\geq 2$.
\end{thm}

\begin{cor} Assume the minimal model program. Let $n>1$. Given $m>0$, there 
exists a Hilbert polynomial such that the moduli space of stable Gorenstein
varieties of dimension $n$ with this Hilbert polynomial has at least $m$ 
components.
\end{cor}

\begin{proof}
Given $m>0$, there exists a positive integer $N$ which factors in at least
$m$ distinct ways as a product of two distinct factors. Choose $m$ pairs 
$(a_i,b_i)$ such that $(a_i-1)(b_i-1)=N$. Let $C_{a_i}$ and $C_{b_i}$ be 
smooth curves of genus $a_i$ and $b_i$, respectively for each $i$. 
 
Let $g_1,\ldots, g_{n-2}$ be distinct integers greater than 1 which are also 
distinct from all of the $a_i$ and $b_i$ and for each 
$j=1,\ldots,n-2$, let $C_{g_j}$ be a smooth curve of genus $g_j$. 
Then the products
\begin{center}
\begin{tabular}{c}
$C_{g_1}\times\cdots\times C_{g_{n-2}}\times C_{a_1}\times C_{b_1}$ \\
 $\vdots$  \\
$C_{g_1}\times\cdots\times C_{g_{n-2}}\times C_{a_m}\times C_{b_m}$
\end{tabular}
\end{center}
have the same numerical invariants, since these can be computed from
the invariants of the curves, and for a product of two curves, $\chi$ and
$K^2$ are both multiples of $(a_i-1)(b_i-1)$.
However, these curves belong to different components of
the moduli space since the genera chosen are distinct.
\end{proof}

One could also draw several easy corollaries of the theorem from the deep 
results in \cite{hm:mc} concerning the moduli spaces of stable curves; in 
particular:

\begin{cor}
$M_{g_1,g_2}$ is of general type if $g_1$ and $g_2$ are distinct and both
greater than 23.
\end{cor}

Also, the rational Picard group is not as simple as in the case
of curves.

\begin{cor}
Let $g_1$ and $g_2$ be distinct integers greater than 2.
Then $\Pic M_{g_1,g_2}\tens \Q\cong (\Pic M_{g_1}\tens\Q)\times(\Pic M_{g_2}
\tens\Q)$.
\end{cor}

\begin{proof}
The moduli spaces of curves are integral schemes of finite type. Furthermore,
$H^1(M_g,\C)=0$ (see, e.g. \cite{ac:modcoho}). Their singularities
are at worst finite quotient singularities, since the Kuranishi spaces for 
curves are smooth and stable curves have a finite automorphism group. Since
finite quotient singularities are DuBois, the results of \cite{dub:rdc} imply
that $H^1(M_g,\C)\rightarrow H^1(M_g,\sh{O}_{M_g})$ is surjective, so the 
latter group is zero. The 
result that the Picard group of the product decomposes as the product of 
Picard groups under these hypotheses is \cite{hart:ag} ex. III.12.6.
\end{proof}

Specifically, for the moduli spaces of curves, the rational Picard group is
freely generated by the Hodge class and the classes of the 
components of the boundary divisor \cite{arb:pic}. The rational Picard group 
of the product has too high a rank for the same to be true. This is 
not surprising, since the cycle structure of the moduli spaces of surfaces is
not as simple as that for curves. In general, the boundary is not likely
a divisor, and there will be other ``geometric'' classes 
which occur, for example, the closure of the locus of surfaces whose canonical
model has rational double points.

\end{document}